\documentclass[11pt]{amsart}

\usepackage{amsmath,amsfonts,amssymb,latexsym,amsthm}
\usepackage{epsfig}
\usepackage{hyperref} %hyperref wants to be loaded last
\usepackage[dvipsnames]{xcolor}

\numberwithin{equation}{section}

\newtheorem{thm}{Theorem}[section]

\newtheorem{conv}[thm]{Convention}

\title{Exploring mazes at random}

\author[Nikita Gladkov, \. Igor Pak]{Nikita Gladkov \. and \. Igor Pak}

 \address{Department of Mathematics, %
	UCLA, Los Angeles, CA 90095, %
	USA}
\email{\{gladkovna,pak\}@math.ucla.edu}

\thanks{\today}

\def\<{\langle}
\def\>{\rangle}

\def\0{{\mathbf 0}}

\def\.{\hskip.06cm}
\def\ts{\hskip.03cm}

\def\nin{\noindent}

\begin{document}

\begin{abstract}
We consider a probabilistic version of the depth-first search on mazes with 
two exits, and show that this algorithm has equal probability of finding 
either exit.  The proof is combinatorial and uses an explicit involution. 
\end{abstract}
\maketitle

% \subsection*{Foreword}
%
Suppose you are entering a maze which has two exits.  You start exploring,
moving from room to room until you find the fiirst exit door.

\smallskip

\nin
{\bf Question:} \. {\em Are you equally likely to find either exit?}

\smallskip

Naturally, the answer depends on what exactly do we mean by a ``maze''
and how exactly do we ``explore'' it.  But most likely, you first instinct
is to say ``No, not necessarily, since one exit might be much closer
to the entrance.''  Hold on to that thought as it would make our
story more interesting.

\subsection*{What is a maze?}
For us, a \emph{maze} is a region in the plane without holes that is aligned
with a square grid, and with some edges representing \emph{walls}.  Without
loss of generality, we can take the region to be an \ts $m \times n$ rectangle.
The \emph{rooms} are the unit squares inside the rectangle.  When two
rooms are not separated by a wall, we say that there is a \emph{door}
between them.

We will always assume that the boundary edges of the maze are
walls except for three: edges $A$, $B$, and $C$ placed counterclockwise
along the boundary.  We also assume that
every room has at least one wall (otherwise it's not really a room, is it?)
Here is a typical example of the kind of maze we consider (only the walls
are drawn, not doors):

	\begin{figure}[hbt]
		\includegraphics[width=13.7cm]{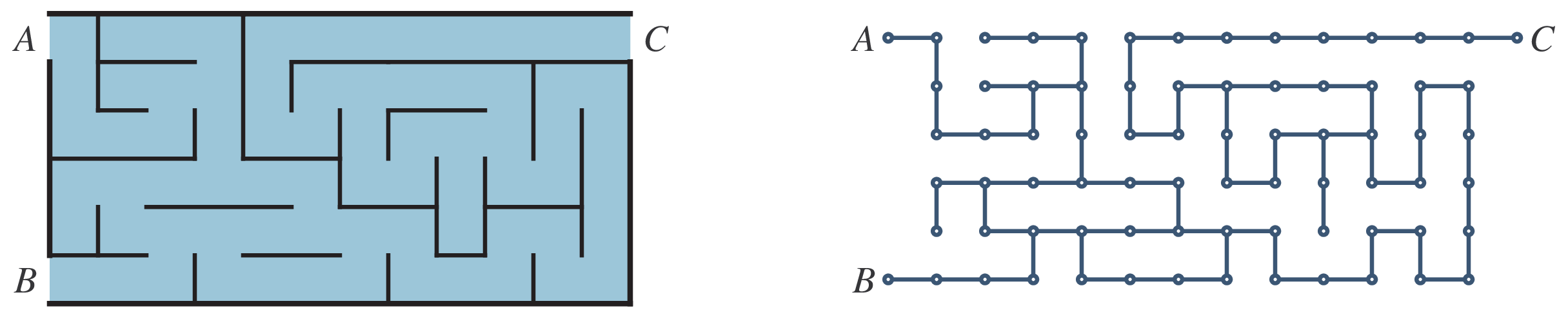}
%		\vskip-.2cm
%		\caption{Maze}
%		\label{f:ex}
	\end{figure}

One can formally define a maze as a planar connected graph with maximal degree 
three, and vertices $A, B,C$  of degree one on the outer boundary, as in the figure.
It will be convenient to keep in mind both formulations.

\subsection*{Shortest paths}
We can now compare shortest paths in this maze, from the entrance at~$A$
to exits at~$B$ and~$C$.  To make this formal, the length of the path is
defined as the
graph distance in the corresponding planar graph.  Comparing the paths in the
figure, it is really hard to imagine that there is an equal probability to
exit at~$B$ and at~$C$:

	\begin{figure}[hbt]
		\includegraphics[width=13.7cm]{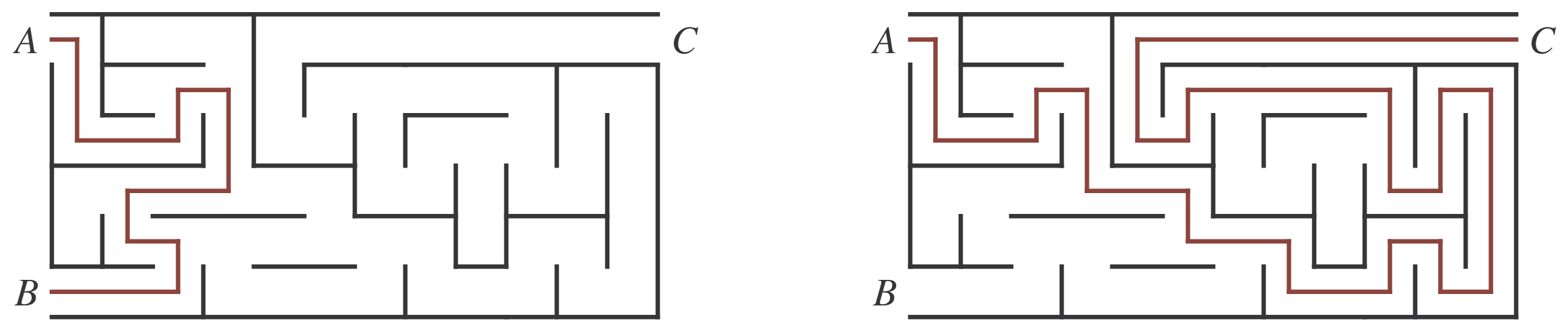}
	\end{figure}

\subsection*{How to explore a maze?}
Let us consider several different ways to explore.

\smallskip

\nin
$(1)$ \. The first approach is deterministic.  Start walking around the maze
in such a way that the right hand always touches the wall.  We call this 
the \emph{right hand on the wall} \ts rule (RHOW).  We can similarly define 
the \emph{left hand on the wall} \ts rule (LHOW).  
In our example, the RHOW gives a path which exits at~$B$, while the 
LHOW gives a path which exits at~$C$.

	\begin{figure}[hbt]
		\includegraphics[width=13.7cm]{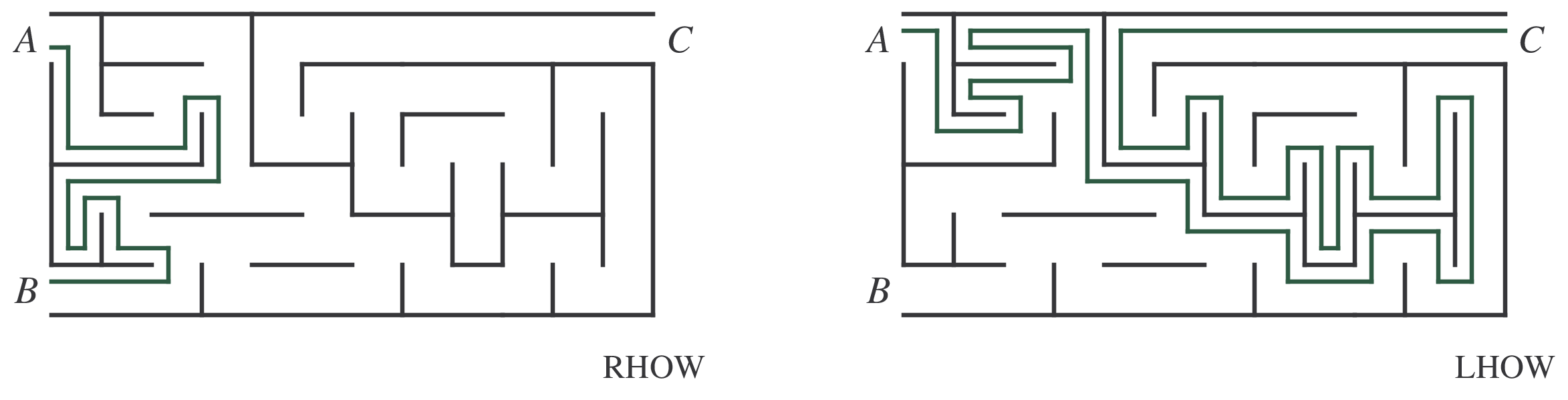}
	\end{figure}

Consider now an \emph{ambidextrous walker} who would need to make a random choice:
should they go with the right hand or left hand?  Say, this choice is random, made by
a fair coin flip.  We suddenly have a probabilistic rule AHOW with equal probability of 
exiting at~$B$ and at~$C$.  Perhaps, this is not overly surprising, but a prelude 
of things to come.

\medskip

\nin
$(2)$ \. Our next approach is algorithmic, and is still deterministic.
Pretend that exits $B$ and $C$ are closed off and explore the whole maze.  
More precisely, walk around the maze in a way to visit every room, 
touch every wall, and walk through every door 
twice --- once in each direction, until returning back to~$A$.  Do this
according to the following local rules which depend on the number of doors
in the room (from~1 to~3) and prior visits:  

\smallskip

	\begin{figure}[hbt]
		\includegraphics[width=9.6cm]{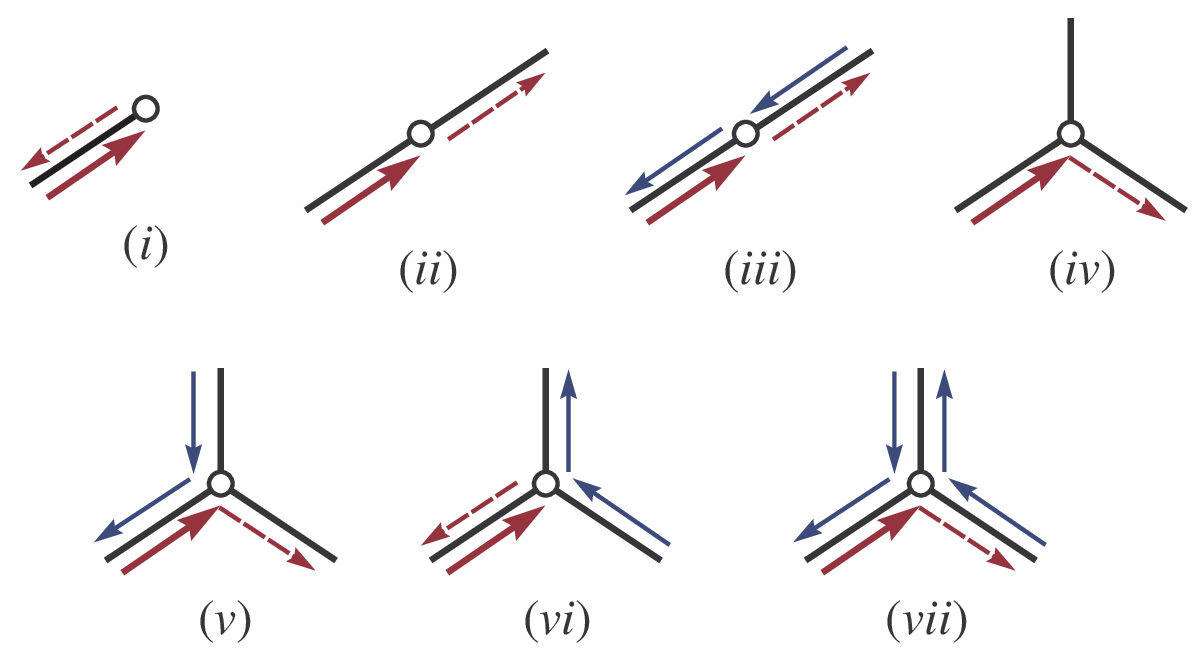}
	\end{figure}

Here the solid red (NE) arrow indicates the door entering a room, dashed red arrow
exits the room, and blue arrow indicate past room visits.  Note that the
only choice here is at rule~$(iv)$, where the rule indicates walking 
through the door on the right from the entry.  
These rules define the \emph{depth-first search} algorithm adapted to
mazes and the right direction in~$(iv)$, see an example in the figure below.  
We denote it RDFS to indicate this preference and the asymmetric nature of 
the algorithm.   

	\begin{figure}[hbt]
		\includegraphics[width=13.7cm]{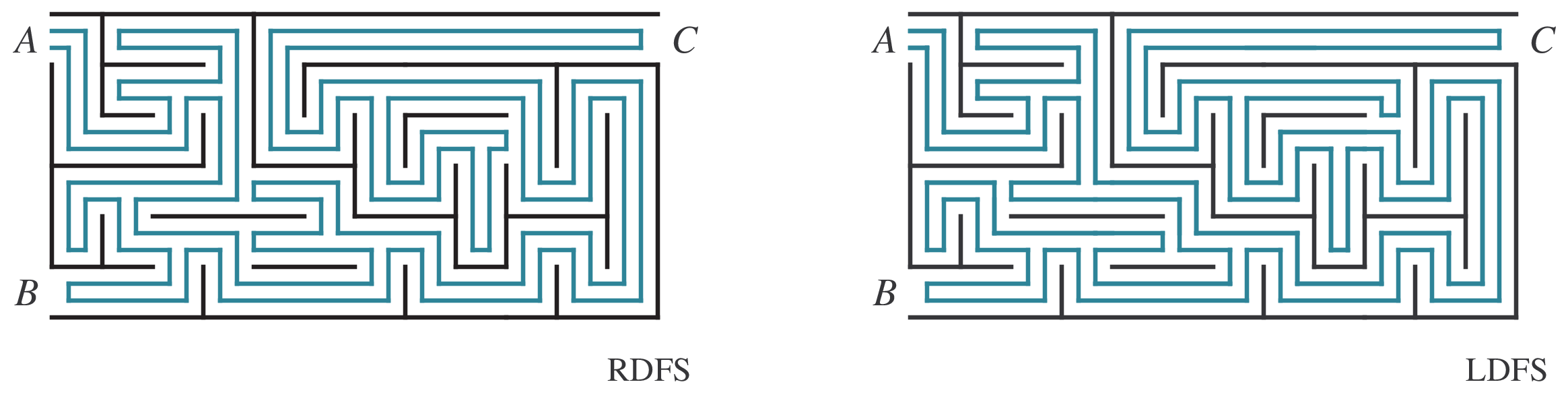}
	\end{figure}

While we explore the maze using RDFS, both exits $B$ and $C$ appear
along the way.  In the example above, we find $B$ before $C$, so we exit
there.
Note that if we replace rule~$(iv)$ with the left direction, 
we get a similar algorithm we denote LDFS.  
In the example, we get $C$ before $B$.  

We can also define the \emph{ambidextrous depth-first search} (ADFS) 
obtained by flipping a fair coin in the beginning and choosing 
one of the two versions.  As with the AHOW, we again have equal 
probabilities of exiting at~$B$ and at~$C$.

\medskip

\nin
$(3)$ \. One can explore a maze via a \emph{random walk}.  At each step,
choose a door in the room uniformly at random and walk there.  This randomized algorithm
will work \emph{eventually}, in a sense that eventually the random walk reaches
either~$B$ or~$C$, but the time that happens can be unbounded.  In this case, the
exit probabilities are not necessarily equal.  The reader might enjoy
calculating the exact probabilities for the example in the figure below.

	\begin{figure}[hbt]
\vskip-.4cm
		\includegraphics[width=4.5cm]{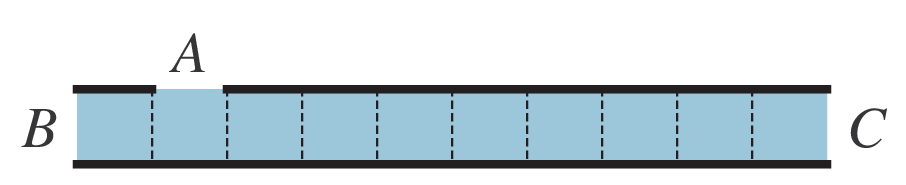}
	\end{figure}

\smallskip

\nin
$(4)$ \. Finally, one can consider \emph{probabilistic depth-first search} (PDFS) 
defined as follows.  The first time the room with three doors is entered, 
flip a fair coin.

$\circ$ \. If \emph{Heads}, use the right door in rule~$(iv)$.

$\circ$ \. If \emph{Tails}, use the left door in rule~$(iv)$.

\nin
This randomized algorithm provably finishes in time bounded by twice the
number of doors in the maze.  For our main example, there are 10 coins being
flipped, so the probability of exiting at $B$ is going to be equal \ts $\frac{a}{2^{10}}$
\ts for some integer~$a$.  The following is our main result:

\medskip

\nin
{\bf Theorem:} \. {\em In {\rm PDFS}, the probabilities of exiting at~$B$ and~$C$ are equal to~$\ts \frac12$\ts. }

\medskip

In other words, when it comes to the ordering of~$B$ and~$C$, 
the PDFS algorithm is similar to the AHOW and ADFS algorithms, 
where only one coin is flipped in the beginning.

\smallskip

\subsection*{The case of a tree maze}
To understand and explain the Theorem, consider the case when all walls
in the maze are connected to the exterior wall.  Note that our main example
is not like that, as it has three connected components of walls not connected
to the exterior wall.  We call these \emph{tree mazes} because the corresponding
graph is a tree, i.e.\ has no cycles.  We use $T$ to denote this tree.

Recall that $A$, $B$ and $C$ are \emph{leaves} (endpoints) of~$T$.  Vertex~$x$ 
in~$T$ of degree three is called a \emph{pivot vertex}, if vertices $A$, $B$ and $C$ 
lie in different connected components when $x$ is deleted. We call these components
\ts $T_A$, \ts $T_B$ \ts and \ts $T_C\ts.$  Observe that every tree~$T$ has a 
unique pivot vertex.  This may seem obvious (see the figure below), 
but the reader might like to carefully prove that.

	\begin{figure}[hbt]
		\includegraphics[width=13.7cm]{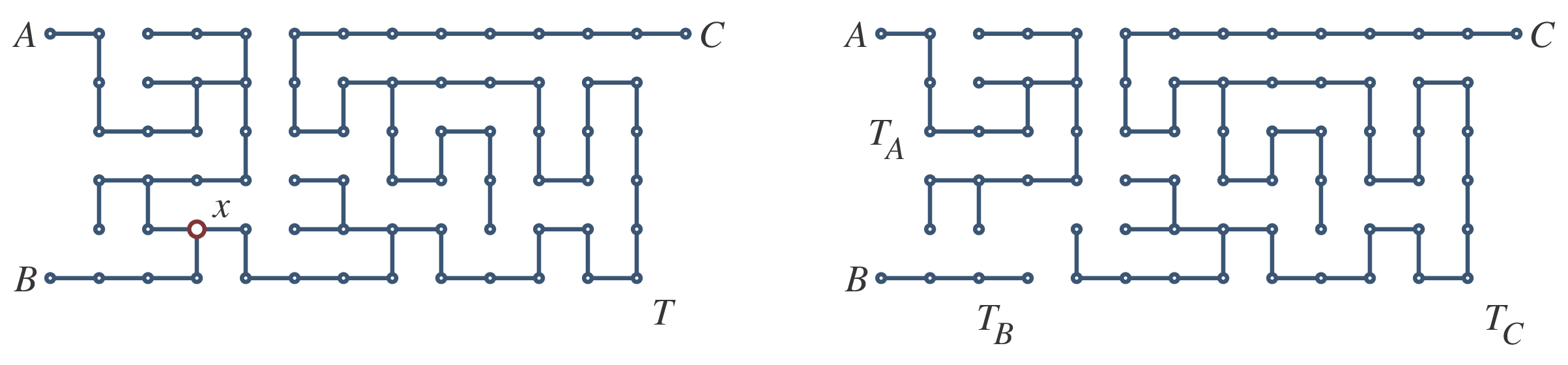}
	\end{figure}

Now, the PDFS starts at $A$ and eventually reaches the pivot vertex~$x$.
Depending on the coin flip, the walk will move to either $T_B$ or~$T_C$ 
with equal probability.  Say, it moved to~$T_B\ts$.  By definition of 
the DFS, we must explore the whole connected component $T_B$ 
before returning back to~$x$.    After that, it fully explores the 
other connected component $T_C\ts$, and only then moves back to~$T_A\ts$.  
In particular, PDFS visits~$B$ \emph{before} it visits~$C$ in this case.  
This is a delicate technical point, the reader might want to think it over.  

We conclude that whether the walk exits at $B$ or at~$C$ depends only
on the outcome of the coin flip at the pivot vertex, which means that these
probabilities are equal.  This proves the Theorem for tree mazes.

\smallskip

\subsection*{General mazes}
We now construct the \emph{DFS tree} as follows.  Every time
a door is opened into a room that was never visited, include
an edge corresponding to that door.  For the RDFS
exploration in our main example, we obtain tree~$T$
as shown in the figure below. 

\smallskip
	\begin{figure}[hbt]
		\includegraphics[width=13.7cm]{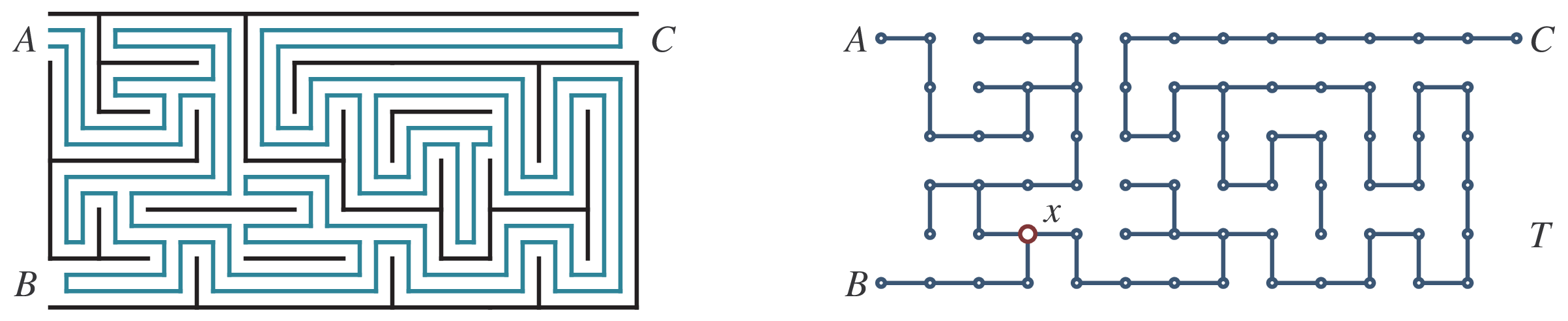}
	\end{figure}

In general, for a PDFS exploration, the DFS tree~$T$ 
depends on the choice of coin flips.  Since rule $(iv)$ is
used only the first time the room is entered and the corresponding
edge is added to~$T$, the PDFS exploration coincides with the 
PDFS exploration on tree~$T$ (for the same outcomes of coin flips).  

From this point on, proceed as in the case of the tree maze.
Find a pivot vertex~$x$ in~$T$.  Change the outcome of the 
coin flip at~$x$.  Observe that the DFS tree corresponding 
to these outcomes of coin flips is also~$T$ (the reader 
might want to explain why).  Therefore, changing the 
coin flip at~$x$ gives an \emph{involution} on all PDFS graph 
explorations that switches explorations where $B$ precedes~$C$, 
with those where $C$ precedes~$B$.  This means that the 
probabilities as in the Theorem are equal again.  This   
completes the proof of the Theorem in full generality. \qed

\smallskip

{\small

\subsection*{Notes on the proof}
First, observe that planarity played no role in the argument.
In fact, the Theorem holds for all connected graphs with maximal degree three.
On the other hand, having maximal degree three and degrees of $A,B,C$
at most one is essential, see exercise~$(1)$ below.

Finally, in the examples above the $A\to B$ portion of the RDFS path 
coincides with the RHOW path.  This is not always the case, 
see exercise~$(2)$ below.
}
	\begin{figure}[hbt]
		\includegraphics[width=13.7cm]{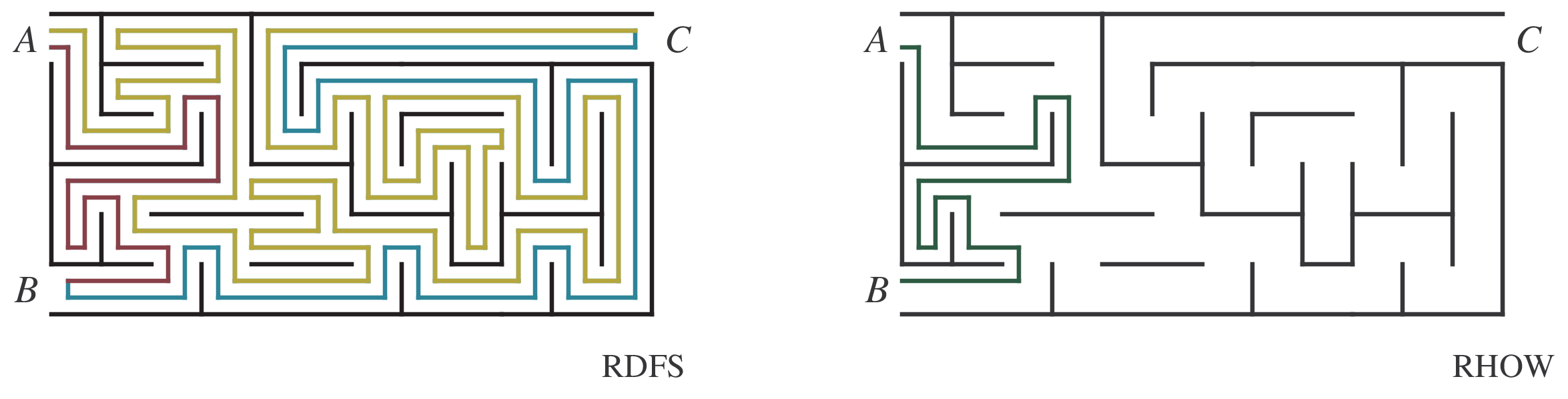}
	\end{figure}

%\smallskip

{\small 

\subsection*{Additional exercises} \.
$(1)$ \. Suppose we have a connected planar graph with maximal 
degree at most three, and vertices $A,B,C$ are on the outer boundary.  
Suppose both $B$ and $C$ are leaves (have degree one), but $A$ has
degree three.  Follow the rules of PDFS, except in the beginning choose
a direction with equal probability~$\frac13\ts.$  
Give an example where the probabilities as in the Theorem are unequal.  
Same question when $A$ and $B$ are leaves, but $C$ has degeee two.

\smallskip

\nin
$(2)$ \. Check that in the example above, the $A\to C$ portion of LDFS 
coincides with the LHOW. Find an example of a maze where this is not true.  
Give necessary and sufficient conditions for these paths to coincide.  

\smallskip

\nin
$(3)$ \. Observe that the RDFS (undirected) path does not always 
coincide with the LDFS path, as in the example above.  Can we have 
both paths with $B$ before~$C$?  

\smallskip

\nin
$(4)$ \. Prove that the RDFS and LDFS can always be drawn in the plane 
without self-intersection.  Are these the only instances of PDFS when
this happens?  

\smallskip

\nin
$(5)$ \. Suppose a maze has three exits: \ts $B$, $C$ and~$D$.
Prove that the probability of exit at $B$ is between \ts $\frac14$ \ts
and \ts $\frac12$\ts.  Does there exist a maze where this probability
is exactly $\frac13$\ts?

\smallskip

\nin
$(6)$ \. In a maze with one entrance $A$ and no exists, consider 
a door connecting two adjacent rooms $x$ and~$y$.  Prove that the 
probability that the door from~$x$ to~$y$ is opened \emph{before} 
the door from $y$ to~$x$, is either \ts $0$, \ts $\frac{1}{2}$ 
\ts or~$\ts 1$.

\smallskip

\nin
$(7)$ \. For a general planar graph, replace every vertex of degree
at least four with a ``ring'' as in the figure below, 
so that the resulting graph has maximal degree three.  
Use the PDFS algorithm for the latter graph to design
the PDFS on the former.  State and prove the analogue 
of the Theorem in this case.

\smallskip

	\begin{figure}[hbt]
		\includegraphics[width=7.8cm]{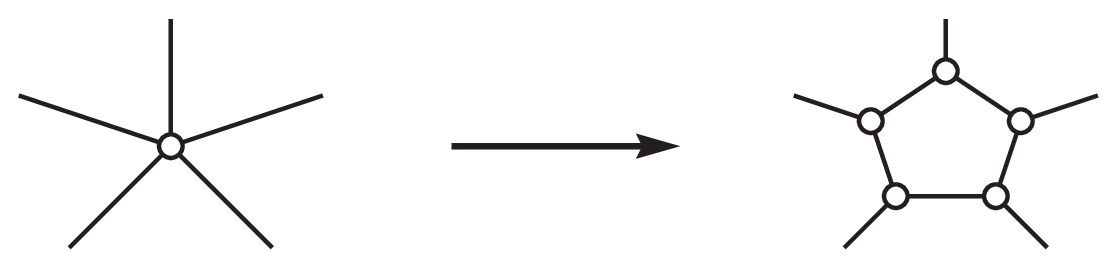}
\vskip-.3cm
	\end{figure}

\subsection*{Historical remarks} \. While solving mazes may seem like a
recreational mathematics, it is the origin of the DFS algorithm, which was
first designed for this purpose (see \cite{Ber}).  We refer to \cite[Ch.~3]{Lucas}, 
for the original discussion of maze exploration attributing the DFS algorithm
to a French mathematician Charles Pierre Tr\'emaux (1859--1882).
See \cite[$\S$22.3]{CLRS} and \cite[Ch.~3]{Even} for a modern treatment
of the DFS and other graph algorithms.
For a random walk on an interval starting at~$A$ with two absorbtion
states $B$ and~$C$, see a careful discussion in \cite[Ch.~XIV]{Feller}.
For mixing times and other aspects of random walks on general graphs, see \cite{LP}.
}

%\newpage

\vskip.7cm

%%%%%%%%%%%%%%%%%%%%%%%%%%%%%%%%%%%%%%%%%%%%%%%%%%%%%%%%%%%%%%%%%%%%%%%%

\end{document}